\documentclass{article}
\usepackage{amssymb}
\newcommand{\eq}{\begin{equation}}
\newcommand{\en}{\end{equation}}

\def\endpf{\hfill $\Box$ \vskip0.5cm}
\def \proof{\noindent{\it Proof.\ }}

\newtheorem{theorem}{\large Theorem} %[section]
\newtheorem{proposition}[theorem]   {\large Proposition}

\newtheorem{corollary}[theorem]{\large  Corollary}

\font\tenmath=msbm10 
\font\sevenmath=msbm7 
\font\fivemath=msbm5 
\newfam\mathfam \textfont\mathfam=\tenmath
\scriptfont\mathfam=\sevenmath \scriptscriptfont\mathfam=\fivemath

\def \\ { \cr }
\def\R{{\mathbb R}}
\def\N{{\mathbb N}}
\def\E{{\mathbb E}}
\def\P{{\mathbb P}}

\def \e{{\rm e}}
\def \f{{\cal F}}

\begin{document}

\title{Asymptotic laws for  nonconservative self-similar fragmentations}

\author{Jean Bertoin
\thanks
{Laboratoire de Probabilit\'es et Mod\`eles Al\'eatoires
et Institut universitaire de France,
Universit\'e Pierre et Marie Curie, 175, rue du Chevaleret,
F-75013 Paris, France}
\hspace{.2cm}
and 
\hspace{.2cm}
Alexander V. Gnedin
\thanks 
{ Mathematisch  Instituut,
Rijksuniversiteit Utrecht,
P.O. Box 80 010, 
3508 TA Utrecht, 
The Netherlands 
}}

\date{\today}
\maketitle

{\bf Abstract}
We consider a self-similar fragmentation process in which the generic particle of size $x$ is replaced at probability rate $x^\alpha$ by its offspring made of smaller particles, where $\alpha$ is some positive parameter.
The total of offspring sizes may be both larger or smaller than $x$ with positive probability.
We show that under certain conditions the typical size in the ensemble is of the order $t^{-1/\alpha}$ and 
that the empirical distribution of sizes converges to a random limit which we characterise in terms of the
reproduction law.

\section{Introduction}

We study the following continuous-time model of particle fragmentation.
Each particle in ensemble is characterised 
by a positive  quantity which we call {\it size}. The generic 
particle of
size $x$  lives a random exponentially distributed time with parameter $x^{\alpha}$,
where $\alpha$ is some fixed real number. 
During the life-time the size does not vary
and at the end the particle splits
into fragments of sizes $x\xi_j$, where
$\{\xi_j\}$ is independent of the lifetime of the particle and follows a
given probability distribution called
  {\it reproduction law}.
Each particle is autonomous, meaning that the splitting probability rate and the 
descendant fragment sizes depend 
only on the size of the particle and not on the history of this or other coexisting particles.
We are interested in the case $\alpha>0$, when particles with smaller size tend to live longer.
We refer to the proceedings \cite{BCP} and the survey \cite{A} for
a number of examples arising in physics and chemistry.

\par The idea of the model was suggested by Kolmogorov in \cite{Kolmogorov}, and the first
results     are due to Filippov \cite{Filippov}. 
Brennan and Durrett \cite{BDI, BDII} 
rediscovered an instance of the model in the context of binary interval splitting.
In a recent series of papers  Bertoin  \cite{Be0, Be1, Be2}
  introduced more involved fragmentation processes in which a 
particle may produce infinitely many generations within an arbitrary  time period
or may undergo a continuous  size erosion; see also \cite{AP, Be, Be3, Miermont, MS, Sch} 
for related examples. 
\par The research so far  
was mainly focussed on the {\it conservative} case  $\sum \xi_j=1$
when the total size is preserved by each splitting. 
It has been shown that the particles demonstrate quite a regular 
long-run behaviour: the typical size in the ensemble  is of the order $t^{-1/\alpha}$ and 
a scaled empirical distribution of sizes converges
to a nonrandom limit.
Filippov \cite{Filippov} proved the convergence of empirical distributions in probability 
(see also \cite{Be2}), while Brennan and Durrett \cite{BDII} showed convergence with
probability one in the binary case. 
Baryshnikov and Gnedin \cite{BG} studied a sequential interval packing problem which may be seen
as a binary instance of {\it dissipative} fragmentation with
 $\sum \xi_j\leq 1$ and ${\mathbb P}(\sum \xi_j< 1)>0$, 
and proved convergence of the mean measures associated with the empirical 
distributions. We mention that some special cases of these mathematical results have also appeared in the literature in physics, see e.g. \cite{BK, KBNG} and references therein.

\par Conservative or dissipative fragmentations can be 
treated both as
continuous-time interval splitting schemes, similar to discrete-time random recursive constructions 
(as in \cite{Mauldin}),
or as state-discretised processes with values in Kingman's partition structures \cite{Be0, Be1, Be2, Be3, BR}.
These approaches fail completely 
if the reproduction law allows  the possibility of size {\it creation},
when the total  size of the offspring may exceed the size of the parent particle. 
Such `improper fragmentations' 
are both physically plausible and useful 
in the situations where the generalised size  
models some nonadditive quantity like, e.g. surface energy by aerosols.

\par By allocating particle of size $x$ at $-\log x$ 
the fragmentation process can be seen  as a branching random walk with location-dependent sojourn times.
From this viewpoint,  a constraint on the sum of sizes seems rather odd, which suggests that 
such a condition is not essential for the asymptotics. 
In this paper we argue that this intuition is indeed correct, in the sense that 
the convergence of properly scaled empirical distributions of sizes holds under fairly general
assumptions on the reproduction law. 
Though we do require that 
the individual offspring sizes 
cannot exceed the parent size, there is no constraint on the total offspring size. 
A new feature appearing in the general nonconservative case 
is that the limit 
of scaled empirical distributions  is not completely deterministic, rather
involves a random factor which admits a characterisation by  a  
distributional fixed-point equation. This phenomenon reminds us, of course,
of strong limit theorems for branching random walks; we shall discuss the connection and
the differences in Section \ref{BRW}.

\par The rest of this work is organised as follows. Notation and basic assumptions are given
in Section \ref{Defin}, and Section \ref{Genea} introduces the genealogical structure and the so-called intrinsic 
martingale  which plays a fundamental role in branching processes. Then we compute the first moment
of power sums and then determine its asymptotics using a contour integral. This
yields the convergence of mean measures in Section \ref{CMM}. An alternative approach based on a limit
theorem of Brennan and Durrett is presented in Section \ref{st-br}. 
The main result of convergence of scaled empirical measures is proved
in Section \ref{L2}. 
Then we provide some examples in Section \ref{examples}, and 
finally, in Section \ref{InFr}, we sketch the extension of the preceding results
to self-similar fragmentations with possibly infinite reproduction measure.

\section{Preliminaries}
\subsection{Definitions and assumptions on the reproduction law}
\label{Defin}
\par The collection $\{\xi_j\}_{j\in\N}$  of offspring sizes of a unit particle
is identified with a random 
sequence of nonnegative real
numbers ranked in the decreasing order (by convention, $\xi_{j}=0$ when $j$ is larger than
the number of children). We shall also view
$\{\xi_j\}_{j\in\N}$ as  a random  set, defined formally as a counting random
measure $\sum\,\delta_{\xi_j}$ on $\,]0,1]\,$. Basically we require
that 
\begin{equation}\label{trivial-ass}
\xi_j\in \,[0,1]\,,\qquad {\mathbb E}\#\{j:\xi_j>0\}>1\,,
\qquad {\mathbb P}(\xi_{1}=1)<1\,.
\end{equation}

\par Many features of the fragmentation process can be expressed in terms of the {\it structural measure}
$$\sigma(B)={\mathbb E}\,\#\{j: \xi_j\in B\}\,,\qquad B\subset \,]0,1],$$
and its Mellin transform
$$\phi(\beta)=\int_0^1 x^{\beta}\sigma({\rm d}x)={\mathbb E}\,\sum_{j=1}^{\infty}\,\xi_j^{\beta}\,$$
which we call the {\it characteristic function}.
In particular, the conditions in (\ref{trivial-ass}) 
amount to the assumptions that 
$\sigma$ is supported by $[0,1]$, that
$\sigma\,]0,1]>1$ and that $\sigma\{1\}<1$.
\par Because
$|\phi(\beta)|\leq \phi({\Re\,}\beta)$,
the natural domain of definition of $\phi$ is 
a complex halfplane to the right of the convergence abscissa $\beta_a$ of the integral.
If $\beta_a=-\infty$ the halfplane is the whole plane, and otherwise the halfplane may be open or closed.
The characteristic function  is analytical in the halfplane,  
strictly decreasing
on the real axis, and in view of (\ref{trivial-ass}) satisfies $\phi(0)> 1$ and $\phi(\beta)\to\sigma\{1\}<1$ 
as $\Re \beta\to\infty$.

\par It is crucial for our results and will be assumed throughout that there exists the {\it Malthusian exponent} 
$\beta^*>0$ satisfying the equation
\begin{equation}\label{char-eqn}
\phi(\beta)=1\,.
\end{equation}
If the Malthusian exponent exists then it is unique and there are no solutions  to
(\ref{char-eqn}) in the halfplane
$\Re
\beta>\beta^*$. And if some $\beta\neq \beta^*$ with $\Re \beta=\beta^*$ satisfies
(\ref{char-eqn}) then $\sigma$  is {\it arithmetic}, meaning that $\sigma$ is a discrete
measure supported by a geometric sequence. Note that in the conservative case $\beta^*=1$ and
in the  dissipative case $\beta^*< 1$.

\par Equation (\ref{char-eqn}) has no real solutions (thus the Malthusian exponent is not defined)   
only if $\phi(\beta_a+)<1$. An example of this situation is  
the measure
$$\sigma({\rm d}x)= c\, {\bf 1}_{\{x<1/2 \}}\,x^{-3/2}\log^{-2} x\,\,{\rm d}x$$
with  a suitable choice of $c$, and $\beta_a=1/2$.
\par Further assumptions about the reproduction law will be introduced in a due place.
Specifically,  the $L^2$-convergence results in Sections \ref{Genea} and \ref{L2} require that
\begin{equation}\label{square}
{\mathbb E}\,\bigg(\sum_{j=1}^{\infty}\, \xi_j^{\beta^*}\bigg)^2<\infty\,.
\end{equation}

\vskip0.5cm

\par {\bf Example.} Consider a dissipative reproduction law induced by the uniform stick-breaking.
Let $U_0,U_1,\ldots$ be i.i.d. uniform, and let
$\{\xi_j\}_{j\in\N}$ stand for the rearrangement in the decreasing order of the sequence 
$$(1-U_j)\prod_{k=0}^{j-1} U_k\,,\qquad j=1,2,\ldots$$
meaning that the uniform portion of the size $1-U_0$ is lost, and the rest is fragmented by the
`random alms' principle (as Halmos called stick-breaking). Equivalently, the offspring $\{\xi_j\}_{j\in\N}$ of a unit particle 
can be seen  as the collection of sizes of a random Poisson-Dirichlet distribution with parameter $1$,
upon removing a size selected by a size-biased pick.
Building the characteristic function
$$\phi(\beta)=\sum_{j=1}^\infty {1\over(1+\beta)^{j+1}}= {1\over \beta(\beta+1)}$$
we see that the abscissa is at  $\beta_a=0$ and the Malthusian exponent is
$\beta^*={(-1+\sqrt{5})/ 2}.$

 \vskip0.5cm

\par It should be noted that there is no substantial constraint 
on $\sigma$ imposed by the requirement that $\sigma$ be a structural measure.
Given  
$\sigma$ on $[0,1]$, satisfying
$\sigma\,]0,1]>1$ and $\sigma\{1\}<1$, a possible reproduction law satisfying (\ref{trivial-ass}) with this structural measure 
can be constructed as follows.
Decompose $\sigma=\sigma_1+\sigma_2$ so that  $\sigma_1$ be probability measure and 
$\sigma_2$ some other measure. Let 
$\xi_1$ be a random point with distribution $\sigma_1$ and let
$\xi_2,\xi_3,\ldots$  be the atoms of a  Poisson point process on the unit interval, with intensity measure
$\sigma_2$. Clearly the point process $\{\xi_j\}$ will have intensity $\sigma$.
\vskip0.5cm

\subsection{Genealogical structure and the intrinsic martingale}
\label{Genea}
In this section, we develop some elements on the genealogical structure of the fragmentation.
 This can be viewed as a different parametrisation of the process,
in which the natural time is replaced by the generation of the different particles.
Specifically, we consider the infinite tree
$${\cal U}\,:=\,\bigcup_{n=0}^{\infty}\N^n\,,$$ with the convention
$\N^0=\{\emptyset\}$. The elements of ${\cal U}$ are called nodes.
 For each $u=(u_1,\ldots,u_n)\in{\cal U}$, we
call $n$ the {\it generation} of $u$ and write $|u|=n$, with the obvious convention
$|\emptyset|=0$. When $n\geq 1$ and
$u=(u_1,\ldots,u_n)$, we call
$u-:=(u_1,\ldots,u_{n-1})$ the father of $u$, and $u,i=(u_1,\ldots,u_n,i)$ the $i$-th child of $u$.

It will be assumed that  $\alpha>0$ and that the fragmentation process $X$ starts from a single particle with unit size, unless explicitly indicated.  
We encode $X$ by putting marks on the nodes of the infinite
tree ${\cal U}$ as follows. The initial particle with unit size corresponds to the ancestor
$\emptyset$, and the mark of $\emptyset$ is the
triple $(\xi_{\emptyset}, g_{\emptyset}, d_{\emptyset})=(1,0,{\bf e})$ where
${\bf e}$ is the instant of the first reproduction, so that $[g_{\emptyset},d_{\emptyset}[$
is the time-interval during which the ancestor particle is alive.
The nodes of the tree at the first generation are used as the labels
of the particles arising at the first split,
i.e. $\xi_1,\ldots,\xi_j, \ldots$. Again, the mark associated to each
of these nodes
$i\in \N^1$, is the triple $(\xi_i, g_{i},d_{i})$, where
$g_{i}=d_{\emptyset}$ stands for the birth-time of the $i$-th child of the ancestor,
and $d_{i}=g_{i}+\xi_{i}^{\alpha}{\bf e}_{i}$ for its death-time. And we iterate the same 
construction with each particle at each generation.

\par
Plainly, the description of the dynamics of fragmentation
entails that its representation as a random marked tree 
enjoys the branching property.
Specifically, 
the distribution of the random marks can be
described recursively as follows:
Given the marks $\left((\xi_{v},g_{v},d_{v}), |v|\leq n\right)$ at nodes of the
first
$n$ generations, the marks at nodes of generation $n+1$ can be expressed in the
form
$$(\xi_{u},g_{u},d_{u})\,=\,(\xi_{u-}\tilde{\xi}_{u},
d_{u-}, d_{u-}+\xi_u^{\alpha}{\bf e}_{u})\,,\qquad u=(u_1,\ldots,u_{n+1}),$$ 
where $u-=(u_1,\ldots,u_n)$ is the father of the node $u$, 
and 

\noindent $\bullet$  when $u-=(u_{1},\ldots,
u_{n})\in\N^{n}$ describes the nodes at the $n$-th generation, 
$(\tilde{\xi}_{u_1,\ldots,u_n,i},i\in\N)$ 
are i.i.d. random sequences distributed according to the reproduction law,

\noindent $\bullet$  when $u$ describes the nodes at the $(n+1)$-th generation, the variables ${\bf e}_{u}$
are i.i.d. standard exponential variables which are independent of the sequences $\tilde{\xi}_{u-,\cdot}$.

\medskip One says that extinction occurs when $\xi_{u}=0$ for all nodes $u$ at generation $n$ for some large enough $n\in\N$. The assumptions  (\ref{trivial-ass}) ensure that
the process $Z_{n}:=\#\{u\in\N^n: \xi_{u}>0\}$ is a supercritical Galton-Watson process (possibly taking the value $\infty$), so the probability of non-extinction is strictly positive.

\medskip
It is well-known from the works of  Jagers \cite{Jag}, Nerman \cite{Ner} and many other authors, that the Malthusian hypothesis for branching processes is connected to
a remarkable martingale. The latter is often referred to as the {\it intrinsic  martingale},
it plays a crucial role in the analysis of the asymptotic behaviour of branching processes.
The following statement, expressed in terms of the genealogical coding, is part of the folklore (we refer to \cite{Roesler} for the last part of the claim). Recall that $\beta^*$ denotes the Malthusian exponent.

\begin{proposition}\label{PSF3} 
Under assumptions {\rm (\ref{trivial-ass})} and {\rm  (\ref{square})}, the process
$$M_n:=\sum_{|u|=n}\xi_u^{\beta^*}\,,\qquad n\in\N$$
is a martingale which is bounded in $L^2$ and in particular, uniformly integrable. 
Its terminal value $M_{\infty}$ is strictly positive conditionally on non-extinction, and satisfies the distributional identity
\begin{equation}\label{FPeqn}
M_{\infty}\stackrel{d}{=} \sum_{j=1}^{\infty}\xi_j^{\beta^*}M^{(j)}_{\infty}
\end{equation}
where $M^{(j)}_{\infty}$ are independent copies of $M_{\infty}$, also independent of $\{\xi_j\}_{j\in\N}$.
The identity taken together with conditions {\rm (\ref{trivial-ass})}, {\rm  (\ref{square})} and ${\mathbb E}M_\infty=1$
 characterises $M_\infty$ uniquely. 
\end{proposition}

\par
{\bf Remark.}
Observe that in the important
case when the reproduction law is conservative, in the sense that 
$\sum_{i=1}^{\infty}\xi_i=1$ a.s., 
we have $\beta^*=1$ and $M_n=1$ for all $n\in \N$, so that the statement is trivial.

\par
\medskip
So far we have described the fragmentation process in terms of its genealogy;
however the problems of interest are often expressed in terms of the natural time scale.
More precisely, the configuration of the fragmentation process at time $t\geq 0$  
consists in the set $\{X_j(t)\}_{j\in\N}$ of particles coexisting at time $t\geq 0$, that is in terms of random point
measures
$$\sum \delta_{X_j(t)}\,=\,\sum_{u\in{\cal U}}{\bf 1}_{\{
g_{u}\leq t < d_{u}\}}
\delta_{\xi_u}\,.$$

Recall we assume that the process starts with a sole  particle of unit size, that is
$X(0)=\{1, 0,\ldots\}$. 
Denoting $X^{(y)}$
the fragmentation process that starts with a  particle of size $y>0$, we have
the fundamental self-similarity identity
\begin{equation}\label{selfsim}
X^{(y)}(t)\stackrel{d}{=}y\,X(ty^{\alpha})\,.
\end{equation}

\par
The intrinsic martingale 
is indexed by the generations of the infinite tree. 
In terms of the time scale of the fragmentation process,
we write $(\f_{t})_{t\geq0}$ for the natural filtration generated by the fragmentation process $X=(X(t), t\geq 0)$.
Proposition \ref{PSF3}
then yields the following result which could also be proved by techniques of so-called stopping lines applied to the branching process induced by the genealogical coding (see e.g. Theorem 6.3 and Corollary 6.6 in Jagers \cite{Jag}).

\begin{corollary}\label{Cint} Under the assumptions of Proposition \ref{PSF3}  for every $t\geq 0$, we have
$$\E\left(M_{\infty}\mid \f_{t}\right)\,=\, M(t,\beta^*):=
\sum_j X_j^{\beta^*}(t)\,,$$
hence the process $ M(\cdot,\beta^*)$ is a square-integrable
$(\f_{t})$-martingale with terminal value $M_{\infty}$.
\end{corollary}

\proof We know that $M_{n}$ converges in $L^2(\P)$ to $M_{\infty}$
as $n$ tends to $\infty$, so
$$\E\left(M_{\infty}\mid \f_{t}\right)=\lim_{n\to\infty}
\E\left(M_{n}\mid \f_{t}\right)\,.$$
On the other hand, it is easy to deduce from the Markov property applied at time $t$
that
$$\E\left(M_{n}\mid \f_{t}\right)=
\sum_{i=1}^{\infty}X_{i}^{\beta^*}(t){\bf 1}_{\{G(X_{i}(t))\leq n\}}
+\sum_{|u|=n}\xi_{u}^{\beta^*}{\bf 1}_{\{d_{u}<t\}}\,,$$
where $G(x)$ stands for the generation of the particle $x$
(i.e. $G(\xi_{u})=|u|$),
and $d_{u}$ for the instant when the particle corresponding to the vertex $u$ splits.
However, for each fixed vertex $u\in{\cal U}$,
$d_{u}$ is bounded from below by the sum of $|u|+1$ independent 
exponential variables which are independent of $\xi_{u}$. It follows that
$$\lim_{n\to\infty }\E \sum_{|u|=n}
\xi_{u}^{\beta^*}{\bf 1}_{\{d_{u}<t\}} =0\,,$$
and we conclude that $\E\left(M_{\infty}\mid \f_{t}\right)=M(t,\beta^*)$.
 \endpf

\section{Asymptotics in mean}

\subsection{Power sums and their means}
\label{PoSum}
 In the sequel, we shall find useful to consider the power-sum functionals
$$M(t,\beta)=\sum_j X_j^{\beta}(t)$$
and their means
$$m(t,\beta)={\mathbb E}M(t,\beta).$$ 
For shorthand we sometimes refer to the size of a particle raised to the power $\beta$ as the {\it $\beta$-size},
thus $M(t,\beta)$ is the total $\beta$-size of the population existing at time $t$.
Two instances with obvious physical interpretations
are the $0$-size equal to the number of particles, and the $1$-size equal to the total
size of the ensemble.

As the set of particles alive at time $t$ is a part of the set of particles that are born
before time $t$, there is the inequality
$$m(t,\beta)=\E \sum_{i=1}^{\infty} X_{i}^{\beta}(t) \leq \E \sum_{u\in{\cal U}}
\xi_{u}^{\beta}{\bf 1}_{\{g_{u}\leq t\}}\,,$$
where $g_{u}$ denotes the birth-time of the particle labelled by the node $u$.
Observe that for each node $u$ at generation $n$, $g(u)$ can be bounded from below
by the sum of $n$ independent standard exponential variables which are also independent
of $\xi_{u}$, and that, in the notation of Section \ref{Defin},
$$\E\sum_{|u|=n}\xi_{u}^{\beta}=\phi (\beta)^n\,,\qquad
\beta>\beta_{a}, n\in\N\,.$$
It now follows from classical large deviation estimates that $m(t,\beta)<\infty$ for every $t\geq 0$ whenever $\beta>\beta_{a}$.

 The first-split decomposition with application of  (\ref{selfsim}) shows that  $M(t,\beta)$ 
satisfies the distributional identity
\begin{equation}\label{d-ident-Mt}
M(t,\beta)\stackrel{d}{=}{\bf 1}_{\{t<d_{\emptyset}\}}+ 
{\bf 1}_{\{t\geq d_{\emptyset}\}}\,\sum_j \xi_j^{\beta}\, M_j(\xi_j^{\alpha}(t-d_{\emptyset}))
\end{equation}
where $d_{\emptyset}$ is the exponential life-time of the progenitor,
 and the $M_j$'s are independent replicas of $M(\cdot,\beta)$,
which are also independent of $d_{\emptyset}$ and $\{\xi_j\}$. 
Computing  expectations we arrive at
the integral equation
$$m(t,\beta)=e^{-t}+\int_0^t e^{-s}\int_0^1 m((t-s)x^{\alpha},\beta)x^{\beta} \,\sigma({\rm d}x)\,.$$  
Differentiating we see that $m(\cdot,\beta)$ is a solution to the Cauchy problem for the integro-differential
equation
\begin{equation}\label{int-eq}
\partial_t \,m(t,\beta)=-m(t,\beta)+\int_0^1 m(x^{\alpha}t,\beta)x^{\beta} \, \sigma({\rm d}x)\,.
\end{equation}
which must be  complemented by the initial value $m(0,\beta)=1\,.$
Uniqueness of $C^{\infty}$ solutions for equations of this type is shown in \cite{Iserles}. 

\par Equation (\ref{int-eq}) defines functions $m$ for all ${\Re}\,\beta>\beta_a$.
For the higher derivatives we have
$$\partial^k_t\, m(t,\beta)=\partial_t^{k} m(0,\beta)\,\,m(t,k\beta +\alpha),$$
thus $m$ is increasing in $t$ for $\beta<\beta^*$ and decreasing for $\beta>\beta^*$.
\par Solving (\ref{int-eq}) in power series is straightforward. 
Introducing
$$\psi(\beta)=1-\phi(\beta)$$
and then
\begin{equation}\label{easygamma}
\gamma(n,\beta)=\prod_{k=0}^{n-1} \psi(\beta+\alpha k)\,
\end{equation}
(by convention, $\gamma(0,\beta)=1$), we compute 
\begin{equation}\label{pow-ser}
m(t,\beta)=\sum_{n=0}^{\infty} {(-t)^n\over n!}\,\gamma(n,\beta)
\end{equation}
which is an entire function of $t\in {\mathbb C}$. It is indeed the right solution because
from the formula for derivatives it is clear that $m(\cdot,\beta)$ should be $C^{\infty}$ for $t\geq 0$. See \cite{Iserles} for yet another  interesting representation of $m(t,\beta)$, in the form of a generalised Dirichlet series.

\subsection{A contour integral}
\label{Integr}
When $\phi(\beta)$ is a rational function, splitting $\gamma(\cdot\,,\,\beta)$ in 
linear factors
shows that the series (\ref{pow-ser}) represents
a generalised hypergeometric function, in which 
case the $|t|\to\infty$ asymptotic expansions have been thoroughly studied by Mellin-Barnes' contour integral 
technique, see \cite{Marichev}. 
This method extends to the more general situation considered here (also see \cite{BG}).

\par Call $\beta$ {\it singular} if $\psi(\beta+\alpha n)=0$ for some integer $n\geq 0$.
For singular $\beta$ the series $m(t,\beta)$ is a polynomial, thus
$m(t,\beta)= \gamma(n,\beta)t^n +O(t^{n-1})$ for $t\to\infty$. 
Analogous  asymptotics hold also for nonsingular $\beta$, but it is more difficult to justify, because
$m$  is then an infinite series  which starts alternating from some term.
A good heuristic amounts to substituting $m\sim c t^a$
into (\ref{int-eq}) -- the left-hand side is then of the order $t^{a-1}$ while the right-hand side is $o(t^a)$ exactly
when $\psi(\beta+a\alpha)=0$, which suggests that $a=(\beta^*-\beta)/\alpha$ is the right exponent.
Although this kind of reasoning can be made precise   it gives no idea of the coefficient, see \cite{Iserles}.

\par Assuming $\beta$ nonsingular we extrapolate the function $\gamma(\cdot\,,\,\beta)$ from the integer values 
to arbitrary complex values $z$ (such that $\phi(\alpha z+\beta)$ is defined) by means of 
the formula
\begin{equation}\label{eqgamma}
\gamma(z,\beta)=\prod_{k=0}^{\infty} {\psi(\beta+\alpha k)\over
\psi(\beta+\alpha(k+z))}\,.
\end{equation} 
Convergence of the product follows as in \cite{BG}, Section 5.  
Thus defined, $\gamma$  satisfies the functional equation
$$\gamma(z+1,\beta)=\psi(\beta+\alpha z)\gamma(z,\beta)$$
 reminiscent of 
the well-known equation for Euler's gamma function.

\par All singularities of the function $\gamma(\cdot\,,\,\beta)$  are the poles located at  roots of (\ref{char-eqn}).
Let 
$${\cal P}_{\beta}=\{z:\,\,\exists\, n \geq 0,\,\, \psi(\beta+\alpha(n+z))=0\}$$
be the set of singular points. 
 Because (\ref{char-eqn}) has no solutions to the right of $\beta^*$, the rightmost point of ${\cal P}_\beta$ is
$z_{\beta}:=(\beta^*-\beta)/\alpha$, where $\gamma(\cdot\,,\beta)$ has a simple pole
provided $\beta^*>\beta_a$.

\par Still assuming $\beta$  nonsingular we have ${\cal P}_{\beta}\cap \{0,1,\ldots\}=\emptyset$.
Since the poles of $\Gamma(-z)$ are nonnegative integers, 
the function
$\Gamma(-z)\gamma(z,\beta)t^z$ (with $t$ as parameter) also has these poles, with residue $(-1)^n t^n \gamma(n,\beta)/n!$
at $z=n$. Defining $\cal C$ to be a vertical line between $\Re z_\beta$ and $n_\beta:=\min(0,\lceil\Re z_{\beta}\rceil)$
we obtain by the residue theorem and an estimate of $\gamma$ 
\begin{equation}\label{oint}
\sum_{n=n_\beta}^\infty {(-t)^n\over n!}\, \gamma(n,\beta)=
{1\over 2\pi {\tt i}}\int_{\cal C} \Gamma(-z)\gamma(z,\beta) z^t\,{\rm d}z\,.
\end{equation}

\par If $\beta^*>\beta_a$ the function $\gamma(\cdot\,,\beta)$ is meromorphic in an open strip containing the line
$\Re z=z_\beta$, and the residue 
at $z_\beta$ is 
$${\rm Res}_{z_\beta}\gamma(z,\beta)=
{\psi(\beta)\over \alpha\psi'(\beta^*)}\,\,
\gamma\left({\beta^*-\beta\over\alpha}\,,\,\alpha+\beta\right)$$
as it follows from the identity 
$$\gamma(s,\beta)={\psi(\beta)\over \psi(\beta+\alpha s)}\,\gamma(s,\alpha+\beta)\,$$
upon expanding the ratio. Replacing $\cal C$ by another integration contour ${\cal C}'$ located
in the half-plane $\Re z< \Re z_{\beta}$ so that all poles of $\gamma(\cdot\,,\,\beta)$ 
in this half-plane lie to the left of ${\cal C}'$
 we obtain the principal-term asymptotics of $m$.
\vskip0.5cm
\begin{theorem}
\label{m-lim} Suppose $\beta^*>\beta_a$, the structural measure $\sigma$ is  nonarithmetic
and conditions {\rm (\ref{trivial-ass})} hold, then
\begin{equation}\label{asymp-m}
m(t,\beta)\sim \Gamma\left({\beta-\beta^*}\over\alpha\right){\psi(\beta)\over \alpha\psi'(\beta^*)}\,\,
\gamma\left({\beta^*-\beta\over\alpha}\,,\,\alpha+\beta\right)~t^{(\beta^*-\beta)/\alpha}\,,\qquad{\rm as~~} t\to\infty.
\end{equation}
for $\Re \beta>\beta_a$.
\end{theorem}
\vskip0.5cm

\par Using the identity
$\gamma(-z\,,\,\alpha z+\beta)\,\gamma(z,\beta)=1$ we can re-write the $\gamma$-factor in (\ref{asymp-m}) as

$$\gamma\left({\beta^*-\beta\over\alpha}\,,\,\alpha+\beta\right)=
{1\over \gamma\left({(\beta-\beta^*)/\alpha}\,,\,\alpha+\beta^*\right)}\,.$$

\par The restriction of $\psi$ to the real segment $]\beta_a,\infty[$ is plainly a concave
increasing function, so the condition $\beta^*>\beta_a$ entails $0<\psi'(\beta^*)<\infty$.
We also remark that in the arithmetic case other poles on the line $\Re z=\beta^*$ would contribute to the coefficient.
Further terms of the asymptotic expansion  can be obtained by pulling the integration contour through other poles
left of $\beta^*$, as long as $\phi$ admits a meromorphic continuation, which can go beyond the convergence abscissa.
 If $\beta^*$ is the only point of $\cal P_{\beta}$ in a closed strip $\beta^*-\theta\leq \Re z\leq \beta^*$
then the rest-term in (\ref{asymp-m}) is estimated as $O(t^{(\beta^*-\beta)/\alpha-\epsilon})$ with
 $\epsilon=\min(1,\theta/\alpha)$.

\subsection{Convergence of the mean measures}
\label{CMM}
We encode the configuration of sizes $X(t)=\{X_j(t)\}$ into the random measure
 $$\sum_j\, X_j^{\beta^*}(t)\delta_{t^{1/\alpha}X_j(t)}.$$
The associated mean measure ${\sigma}^*_t$ is defined by  the formula
\begin{equation}\label{bound-f}
\int_0^{\infty} f(x){\sigma}^*_t({\rm d}x)={\mathbb E}\sum_j f(t^{1/\alpha}X_j(t)) X_j^{\beta^*}(t)
\end{equation}
which is required to hold for all compactly supported continuous functions $f$. 
It is easily seen that ${\sigma}^*_t$ is a probability measure. 
Our next goal is to show that the measures ${\sigma}^*_t$ converge weakly to a probability measure 
$\rho$ on $]0,\infty[\,$.

\par Because
$$t^{(\beta-\beta^*)/\alpha}m(t,\beta)=\int_0^\infty x^{\beta-\beta^*}{\sigma}^*_t({\rm d}x)$$
the convergence of $t^{(\beta-\beta^*)/\alpha}m(t,\beta)$ implied by (\ref{asymp-m})
amounts to the convergence of power moments
$$\int_0^\infty x^{\beta-\beta^*}{\sigma}^*_t({\rm d}x)\to \Gamma\left({\beta-\beta^*}\over\alpha\right){\psi(\beta)\over \alpha\psi'(\beta^*)}\,\,
{1\over \gamma\left({(\beta-\beta^*)/\alpha}\,,\,\alpha+\beta^*\right)}\,.~$$
Specialising   $\beta=\beta^* +\alpha k$ this simplifies to 
$$
\int_0^1 x^{k\alpha}\,{\sigma}^*_t({\rm d}x)\to {(k-1)!\over \alpha\psi'(\beta^*)}
\prod_{j=1}^{k-1}{1\over \psi(\beta^*+\alpha j)}
$$
in application of (\ref{eqgamma}). It is easy to check and is well known \cite{BY, BDII}
that for $n=0,1,\ldots$  the related moment problem is determinate, whence the following
result.
\vskip0.5cm

\begin{theorem}
\label{weak-conv} Under assumptions of {\rm \,\,Theorem \ref{m-lim}}
the measures ${\sigma}^*_t$ converge weakly, as $t\to\infty$, to a probability 
measure $\rho$ uniquely determined by its
power moments
\begin{equation}\label{moments-rho}
\int_0^\infty x^{k\alpha}\rho({\rm d}x)={(k-1)!\over \alpha\psi'(\beta^*)}
\prod_{j=1}^{k-1}{1\over \psi(\beta^*+\alpha j)}\,,\qquad  k=1, 2, \ldots
\end{equation}

\end{theorem}

\vskip0.5cm

\par Theorems \ref{m-lim} and \ref{weak-conv} imply that 
\begin{equation}\label{pow-conv}
t^{(\beta-\beta^*)/\alpha} m(t,\beta)\to\int_0^\infty x^{\beta-\beta^*}\,\rho({\rm d}x)\,
\end{equation}
which extends the convergence of expectations in (\ref{bound-f}) to a wider class of functions $f$.

\subsection{Self-similar stick-breaking process}
\label{st-br}
A key tool in the conservative case treated in \cite{BDI, BDII, Be0, Be1, Be2} has been the following observation 
related, somewhat paradoxically, to the simplest dissipative case of a singleton ensemble 
with reproduction law $\{\eta\}$,
where $\eta$ is a  random variable assuming values in $\,]0,1]$.
That is to say, if at some time the particle has size $x$ then, independently of the history,
the particle shrinks with probability rate $x^{\alpha}$
and the new size after the shrink becomes $x\eta$ where $\eta$ follows $\widehat{\sigma}$,
and  $\widehat{\sigma}$ is the probability 
law of $\eta$.

\par We recollect briefly a result from \cite{BDII} (also see \cite{Be2, BC}).  Introduce 
$$\widehat{\psi}(\beta)=1-\int_0^1 x^{\beta}\,\widehat{\sigma}({\rm d}x)\,$$
and suppose $\widehat{\psi}(0+)<\infty$.
 Let $L_t$ be the sole size at time $t$, and $\widehat{m}(t,\beta)={\mathbb E}L_t^\beta$.
Assuming $\widehat{\sigma}$ non-arithmetic, Brennan and Durrett \cite{BDII} proved that for $t\to\infty$
$$t^{1/\alpha}L_t\stackrel{d}{\to} Y^{1/\alpha}\, \qquad{\rm and~~~~~} \widehat{m}(t,\beta)\to {\mathbb E}\, Y^{\beta/\alpha}$$
where $Y$ is a random variable with moments
$${\mathbb E}\,Y^k={(k-1)!\over \alpha \widehat{\psi}'(0+)} \prod_{j=1}^{k-1}{1\over \widehat{\psi}(\alpha j)}.$$
The convergence of moments $\widehat{m}(t,\beta)$ was shown for all real $\beta$ strictly to 
the right of the convergence abscissa of $\widehat{\psi}$ (which is nonpositive due to the normalisation
$\widehat{\sigma}\,]0,1]=1$).
They also suggested the explicit  representation 
\begin{equation}\label{Y-repr}
Y\stackrel{d}{=}\sum_{k=0}^\infty \epsilon_k\prod_{j=0}^k \eta_j^\alpha
\end{equation}
where all $\eta_j, \epsilon_j$ are independent, $\epsilon_j$ is mean one exponential, $\eta_j$ for $j>0$ are replicas of 
$\eta$, and $\eta_0$ follows the law
$${\mathbb P}(\eta_0\in {\rm d}x)={\widehat{\sigma}[x,1] {\rm d}x\over \widehat{\psi}'(0) x}\,\qquad x\in \,]0,1].$$
Each  product in (\ref{Y-repr}) corresponds to the size of the particle after $k$ splits, conditionally given the
initial size is $\eta_0$, thus
formula (\ref{Y-repr}) identifies $Y$ with the well-known exponential functional of a (stationary) 
compound Poisson process, see e.g. the survey
\cite{BY}. 

\vskip0.5cm
\par In the conservative fragmentation case, the above `stick-breaking' process describes
the evolution of a particle tagged by an atom of isotope that was injected at a random  uniform
location  into the progenitor unit size.
The mechanism which determines the line of descent of the tagged particle amounts, at each consecutive split, 
to a random size-biased
pick from the child particles. Thus,  defining
$\widehat{\sigma}({\rm d}x):=x \,\sigma({\rm d}x)$  to be the distribution 
of a size-biased pick from $\{\xi_j\}$, 
we obtain the relation $\widehat{m}(t,\beta-1)=m(t,\beta)$ which was observed in   \cite{BDII}, p. 112.
\par In the general nonconservative case choosing
$\widehat{\sigma}({\rm d}x):=x^{\beta^*}\sigma({\rm d}x),$
we still get 
$\widehat{m}(t,\beta-\beta^*)=m(t,\beta)$ and $\widehat{\psi}(z)=\psi(z+\beta^*)$.
An interpretation of these relations akin to the tagged fragment process
is available through the so-called {\it spine} which appears in
the conceptual approach of Lyons {\it et al.}
\cite{LPP, Ly} for convergence of martingales in branching processes. 
The measure $\rho$ may be identified with the distribution of $Y^{1/\alpha}$. 
A consequence of this discussion and the result of Brennan and Durrett is the following corollary.
\vskip0.5cm
\begin{corollary} The conclusion of {\rm Theorem \ref{weak-conv}} remains valid even if 
the assumption $\beta^*>\beta_a$ is replaced
by the weaker $\psi'(\beta^*+)<\infty$ 
\end{corollary}
\vskip0.5cm

\par The tiny improvement upon Theorem \ref{weak-conv} 
appears in the case
where the characteristic function is defined in a closed half-plane
 and 
$\beta^*=\beta_a$, i.e. the Malthusian exponent falls exactly on the convergence abscissa.
An example of such situation is the structural measure of the form 
$\sigma({\rm d}x)=c {\bf 1}_{\{x<1/2\}} x^{-3/2} |\log x|^{-3}\,{\rm d}x\,$ with a suitable
$c$. The method based on contour integration  requires in such cases the
analytical continuation of $\phi$ in a domain to the left of $\beta_a$.

\par Alternatively, along the lines in
\cite{BDI, BDII, BC}, the
renewal theory can be applied also in the nonconservative case, to  prove first  the convergence of measures
$\sigma_t^*\to \rho$ 
and then to justify the asymptotics of mean power-sums using the uniform integrability.  

\section{Limit theorem for the empirical measure of the fragments}
\subsection{$L^2$-convergence}
\label{L2}
Our principal result improves on the convergence of the mean measures $\sigma_t\to\rho$ in Theorem \ref{weak-conv}  and says that the 
scaled empirical measures induced by $X(t)$ converge in a $L^2$-sense to the measure $M_\infty \rho$, where $M_{\infty}$ is the terminal value of the intrinsic martingale
(cf. Proposition \ref{PSF3}).

\vskip0.5cm
\begin{theorem}
\label{strongconv} Assume {\rm (\ref{trivial-ass})}, {\rm  (\ref{square})}, that
$\beta^*>\beta_a$ and that 
$\sigma$ is nonarithmetic. Then
for any bounded continuous $f$
$$
{\rm {\it L}^2\!-\!}\lim_{t\to\infty} \sum_j X_j^{\beta^*}(t)f(t^{1/\alpha}X_j(t))=M_\infty\,\int_0^{\infty}f(x)\rho({\rm d}x).
$$
\end{theorem}

\proof We need to show that 
\begin{equation}\label{rely-on-Jean}
{\mathbb E}\left(\sum_{i,j}  X_i^{\beta^*}(t)f(t^{1/\alpha}X_i(t))
X_j^{\beta^*}(t) g(t^{1/\alpha}X_j(t))
\right)\to 
{\mathbb E}M_\infty^2\left(\int_0^\infty f(x)\rho({\rm d}x)\right)\left(\int_0^\infty g(x)\rho({\rm d}x)\right)
\end{equation}
for positive $f$ and $g$  bounded from above by $1$. 
Indeed, suppose (\ref{rely-on-Jean}) is shown. Denote
$$A_t=\sum_j X_j^{\beta^*}(t) f(t^{1/\alpha}X_j(t)).$$
Take $f=g$ to conclude from (\ref{rely-on-Jean}) that
$$\lim_{t\to\infty} {\mathbb E} A^2_t={\mathbb E} M^2_\infty \left(\int_0^{\infty} f(x)\rho({\rm d}x)\right)^2.$$
Similarly, setting $g=1$
$$\lim_{t\to\infty} {\mathbb E}(A_tM(t,\beta^*))={\mathbb E}M^2_\infty \int_0^{\infty} f(x)\rho({\rm d}x).$$
Recalling from Corollary \ref{Cint} that ${\mathbb E}M(t,\beta^*)^2\to {\mathbb E}M^2_\infty$ and combining the above we get the desired 
$$\lim_{t\to\infty}{\mathbb E}\left(A_t-M(t,\beta^*)\int_0^\infty f(x)\rho({\rm d}x)\right)^2=0\,.$$

\par To prove (\ref{rely-on-Jean})
let us replace $t$ by $t+s$ and condition on the configuration of sizes $X(s)$.
At time $t+s$ two coexisting particles may stem from the same ancestor that lived at time $s$ 
or from two different ancestors;
write $i\sim_s j$ in the first case, and write $i\not\sim_s j$ in the second.
The  sum  in the left-hand side of (\ref{rely-on-Jean}) is split then in two
$$S_1+S_2={\mathbb E}\left(\sum_{i\sim_s j}\cdots|X(s)\right)+
{\mathbb E}\left(\sum_{i\not\sim_s j}\cdots|X(s)\right).$$
Using the fundamental self-similarity relation (\ref{selfsim}) and the Markov nature of the fragmentation process
we estimate the first sum as 
$$S_1\leq\sum_k X_k^{2\beta^*}(s)~{\mathbb E}\left(\sum_j X_j^{\beta^*}(t)\right)^2$$
hence by (\ref{asymp-m}) and Corollary \ref{Cint}
$${\mathbb E}S_1<{\rm const}\,s^{-\beta^*/\alpha}\to 0\qquad {\rm as~~}s\to\infty$$
uniformly in $t$.

\par Dealing with $S_2$ requires more effort. We use the parallel notation $y_j=X_j(s)$.
Write $i\searrow k$ if the size $X_i(t+s)$ stems from $y_k$.
By independence, the descendants of different  particles with sizes $y_k$ and $y_{\ell}$ evolve independently, thus
grouping the  sizes $X_j(t+s)$ by the ancestors at time $s$ yields
$$S_2=\sum_{k\neq \ell}\left({\mathbb E}\sum_{i\searrow k}\cdots\right)
                       \left({\mathbb E}\sum_{j\searrow \ell}\cdots\right)$$
However, by self-similarity and convergence of the mean measures
$${\mathbb E}\sum_{i\searrow k} y_k^{\beta^*} X_i^{\beta^*}(ty_k^{\alpha})f((t+s)^{1/\alpha}y_kX_i(ty_k^{\alpha}))
\to y_k^{\beta^*}\int_0^\infty f(x)\rho({\rm d}x)$$
$${\mathbb E}\sum_{j\searrow \ell} y_\ell^{\beta^*} X_j^{\beta^*}(ty_\ell^{\alpha})g((t+s)^{1/\alpha}y_\ell X_j(ty_\ell^{\alpha}))
\to y_\ell^{\beta^*}\int_0^\infty g(x)\rho({\rm d}x)$$
as $t\to\infty$, therefore by dominated convergence

$${\mathbb E} S_2\sim \left(\int_0^\infty f(x)\rho({\rm d}x)\right) \left(\int_0^\infty g(x)\rho({\rm d}x)\right) 
{\mathbb E}\sum_{k\neq\ell} X_k^{\beta^*}(s)X_{\ell}^{\beta^*}(s)$$
as $s\to\infty$.
It remains to note that 
$${\mathbb E}\sum_{k\neq\ell} X_k^{\beta^*}(s)X_{\ell}^{\beta^*}(s)
\sim {\mathbb E}\sum_{k\,,\,\ell} X_k^{\beta^*}(s)X_{\ell}^{\beta^*}(s)={\mathbb E}M^2_s\to {\mathbb E}M^2_\infty$$
because
$${\mathbb E}\sum_{k} X_k^{2\beta^*}(s)=m(t,2\beta^*)\to 0\,.$$
\endpf

\vskip0.5cm

\par 
{\bf Remarks.} We mention that if we replace the assumption (\ref{square}) by the weaker
$${\mathbb E}\,\bigg(\sum_j\, \xi_j^{\beta^*}\bigg)^p<\infty\,$$
for some $1<p\leq 2$, the calculation of the expectation of the
sum of the $p$-th powers of  jumps then shows that the intrinsic martingale
$M_n$ in Proposition \ref{PSF3} is  bounded in $L^p$, see e.g. Neveu \cite{Nev}. Then techniques of Nerman \cite{Ner}, based on generalisations of the law of large numbers for branching processes, can be adapted to extend Theorem \ref{strongconv} to this situation (of course, $L^2$ convergence then has to be replaced by $L^p$-convergence).

\par Finally, it is known that in the binary conservative case the scaled  empirical measures
converge with probability one
\cite{BDII}. It would be interesting to extend this result to the nonconservative case.

\vskip0.5cm
\par In parallel to (\ref{pow-conv}) there is the following extension of Theorem \ref{strongconv} to power functions.
\vskip0.5cm

\begin{corollary} Under assumptions of {\rm Theorem \ref{strongconv}}
$$L^2\!-\!\lim_{t\to\infty} t^{(\beta-\beta^*)/\alpha} \sum_j X_j^{\beta}(t)= M_\infty \int_0^\infty x^{\beta-\beta^*}\rho({\rm d}x)$$
for $\Re\, \beta>\beta_a$.
\end{corollary}

\proof Along the same line, the proof is reduced to showing that
$$\sup_{t\geq 0} 
t^{(\beta^*-\beta)/\alpha} \,y^{-\beta^*}\,{\mathbb E}\sum_j \left(X_j^{(y)}(t)\right)^{\beta}$$
is bounded uniformly in $y\in \,]0,1[\,$.
And the latter follows by noting that 
(\ref{selfsim}) and (\ref{pow-conv}) imply 
$$ t^{(\beta-\beta^*)/\alpha} \,{\mathbb E}\sum_j \left(X_j^{(y)}(t)\right)^{\beta}\to y^{\beta^*}\int_0^\infty
x^{\beta-\beta^*}\,\rho({\rm d}x)\,.$$
\endpf

\subsection{Comparison with Branching Random Walks in continuous time}\label{BRW}
In this section, we compare the case of positive indices of self-similarity $\alpha>0$ which we have considered so far, with the simpler {\it homogeneous} case when $\alpha=0$.

So suppose here that $\alpha=0$. The process $Z(t)=\{Z_j(t)\}$ with $Z_j(t)=-\log X_j(t)$ is then a continuous-time {\it branching random walk}, as studied 
in \cite{Uch, Bi2}. We have in the notation of Section \ref{PoSum} that $\gamma(n, \beta)= \psi^n(\beta)$ and the identity (\ref{pow-ser}) for the moments of power sums thus becomes $m(t,\beta)=\exp(-t\psi(\beta))$, a formula which is well-known in the context of branching random walk; see  \cite{Uch, Bi2}.

It is read from the work of Biggins \cite{Bi1, Bi2} that for
every $\beta>\beta_a$, the process
$$\,W(t,\beta)\,:=\,\exp(t\psi(\beta))\sum_{j}X_j^{\beta}(t)\,,\qquad
t\geq0$$ is a martingale with c\` adl\` ag paths, which
converges almost surely and in their mean as
$t\to\infty$ provided that $\psi(\beta)<\beta\psi'(\beta)$, see \cite{BR} for details.
Note that this holds in the special case when $\beta=\beta^*$ is the Malthusian exponent, for which there is the identity $W(t,\beta^*)=M(t,\beta^*)=M_{t}$.

\par When furthermore the reproduction law is not arithmetic, 
the asymptotic behaviour of the empirical distribution of sizes can be described as follows: 
for every continuous function
$f:\R\to\R$ with compact support,
\begin{equation}\label{eqbiggins}
\lim_{t\to\infty}  \sqrt t \, 
\e^{-t(\beta\psi'(\beta)-\psi(\beta))}\, \sum_j\,f(t
 \psi'(\beta)-Z_j(t))\,=\, {W(\infty,\beta)\over \sqrt{2\pi
|\psi''(\beta)|}}\int_{-\infty}^{\infty}f(-z)\e^{\beta z}{\rm d}z\,.
\end{equation}
where $W(\infty,\beta)$ is the terminal value of $W(t,\beta)$. See \cite{Uch, Bi2}.

\medskip \par
Theorem \ref{strongconv} bears obvious similarities with (\ref{eqbiggins}). It is interesting to observe that in the homogeneous case $\alpha=0$, sizes decay exponentially fast and the limiting scaled empirical measure is always exponential (up-to a random factor),  whereas for $\alpha>0$ the decay of sizes is
power-like and the limiting scaled empirical measure depends crucially on the structural
measure $\sigma$ (more precisely, $\sigma$ can be recovered from the limiting scaled empirical measure for $\alpha>0$, but not for $\alpha=0$).

\section{Examples} \label{examples}

\subsection{Filippov's example revisited}
\label{Filippov}
Extending an example in Filippov (see \cite{Filippov}, section 8) consider the structural measure 
$$\sigma({\rm d}x)=\lambda x^{\theta-1}\,{\rm d}x\,,\qquad x\in \,]0,1]$$
with parameters $\lambda>\min(\theta, 0)$ and {\it arbitrary} $\theta\in {\R}$.
For $\lambda<\theta+1$ the fragmentation is mean-value dissipative. We have 
$$\phi(\beta)={\lambda\over \theta+\beta}\,,\qquad \beta^*=\lambda-\theta\,, \qquad \psi(\beta)= {\beta-\beta^*\over\beta+\theta}.$$
The characteristic function is thus meromorphic in $\mathbb C$ with a unique simple pole at $\beta_a=-\theta$.

\par Computing 
$$\gamma(n,\beta)={(A)_n\over (B)_n}$$
we see that this is the ratio of two Pochhammer factorials, with $A=(\beta-\beta^*)/\alpha\,$ and $B=(\theta+\beta)/\alpha$,
thus $m(t,\beta)=~  _1F_1(A;B;-t)$ is Kummer's hypergeometric function. The analytical extension of $\gamma$ is
$$\gamma(z,\beta)={\Gamma(A+z)\Gamma(B)\over \Gamma(B+z)\Gamma(A)}.$$

\par Computing the moments we obtain 
$$\int_0^\infty x^{\alpha k}\,\rho({\rm d}x)=(\lambda/\alpha)_k$$
which identifies $\rho$ as
$$\rho({\rm d}x)={\alpha\over \Gamma(\lambda/\alpha)}\, x^{\lambda-1}e^{-x^\alpha}\,{\rm d}x\,,\qquad x \geq 0\,.$$
Note that  the shape parameter $\theta$ cancels and does not appear in $\rho$. It follows that
$\sigma_t:=x^{-\beta^*}\sigma^*_t$ (the intensity of $M \,\delta_{t^{1/\alpha}X_j(t)}$) converges to the measure
 $$x^{-\beta^*}\rho({\rm d}x)=\,{\alpha\over \Gamma(\lambda/\alpha)}\, x^{\theta-1}e^{-x^\alpha}\,{\rm d}x,\qquad x \geq 0\,$$
in accord with the case $\lambda=2,\, \theta=1$ considered in \cite{BDII} in connection with
the conservative binary fragmentation with $\xi_1$ uniform and $\xi_2=1-\xi_1.$
\par It follows that the mean number of particles satisfies
$$m(t,0)\sim {\Gamma(\theta/\alpha)\over \Gamma(\lambda/\alpha)}\,\,t^{(\lambda-\theta)/\alpha}\,,\qquad t\to\infty\,.$$
which agrees with a special case in \cite{Filippov}.
Of course, this formula makes sense only for $\theta>0$, because our asymptotics for $m(t,\beta)$ hold only for
$\Re \beta>\beta_a$, thus for $\theta\leq 0$ the value $\beta=0$  is not considered.

\vskip1cm
 
 \par The case $\lambda=1, \,\theta=0$, when $\sigma({\rm d}x)=x^{-1}{\rm d}x$ corresponds 
to the conservative fragmentation generated by the uniform stick-breaking,
as in the example in Section \ref{Defin} (but without removing a piece).
It is well known that the distribution of a size-biased pick from $\{\xi_j\}$ is uniform, and this implies
that the intensity measure of $\sum \delta_{\xi_j}$ is indeed $\sigma({\rm d}x)=x^{-1}{\rm d}x$.

\subsection{Hypergeometrics}
\label{Hyperg}
The following is a further generalisation of Filippov's example, and covers the class of dissipative binary fragmentations
treated in \cite{BG}.
 Consider a Dirichlet polynomial
\begin{equation}\label{qp}
g(x)=\sum_{j=1}^p \lambda_j \,x^{\theta_j-1}
\end{equation}
which is non-negative on $\,]0,1]$ and has 
real parameters satisfying 
$$\sum_{j=1}^p {\lambda_j\over \theta_j}>1\,.$$ 
Then $\sigma({\rm d}x)=g(x){\rm d}x$ is a measure on $\,]0,1]$ with rational
characteristic function 
$$\phi(\beta)=\sum_{j=1}^p {\lambda_j\over \theta_j+\beta}$$
and by the assumption the rightmost root of $\phi(\beta)=1$ is positive, denote it also $\beta_1=\beta^*$
and denote further roots $\beta_2,\ldots, \beta_p$ (the roots are certainly different from the poles of $\phi$).
\par Observe that
$$\psi(\beta)=\prod_{j=1}^p {\beta-\beta_j\over \beta+\theta_j}\,,$$
thus assuming $\alpha=1$ (without loss of generality) we have
$$m(t,\beta)=\sum_{n=0}^\infty {(-t)^n\over n!} \prod_{j=1}^p {(\beta-\beta_j)_n\over (\beta+\theta_j)_n}$$
where we recognise a generalised hypergeometric function of the type $_p F_p$.
By (\ref{m-lim}) we have $m(t,\beta)\sim c(\beta)t^{\beta-\beta^*}$ for $\Re\beta>\beta_a$.
Noting that

$$\psi'(\beta^*)={1\over\beta^*+\theta_1}\prod_{j=2}^p {\beta^*-\beta_j\over \beta^*+\theta_j}$$
and manipulating infinite products, the
coefficient is evaluated in terms of the gamma function as 
$$c(\beta)=\prod_{j=2}^p {\Gamma(\beta^*-\beta_j)\over \Gamma(\beta-\beta_j)}
\prod_{j=1}^p{\Gamma(\beta+\theta_j)\over\Gamma(\beta^*+\theta_j)}\,.$$ 
This allows to recover 
the density by Mellin inversion as
$${{\rm d}\rho\over{\rm d}x}={1\over 2\pi{\tt i}}
{\prod_{j=2}^p \Gamma(\beta^*-\beta_j)\over\prod_{j=1}^p \Gamma(\beta^* + \theta_j)}
\int_{-{\tt i}\infty}^{{\tt i}\infty}
{\prod_{j=1}^p \Gamma(z+\beta^*+\theta_j)\over \prod_{j=2}^p \Gamma(z+\beta^*-\beta_j)}\,\,x^{-z-1}\,{\rm d}z\,,\qquad x\geq 0$$
which is an instance of Meijer's $G$-function, see \cite{Marichev}.

\par The limit measure is uniquely determined by the integer moments which can be written as
$$\int_0^\infty x^k\,\rho({\rm d}x)={(k-1)!\over \psi'(\beta^*)}
\prod_{j=1}^p
{(\beta^*+1+\theta_j)_{k-1}\over (\beta^*+1-\beta_j)_{k-1}}\,,$$
where  the derivative may be also computed as
$$\psi'(\beta^*)=\sum_{j=1}^p {\lambda_j\over (\beta^* +\theta_j)^{2}}\,.$$

\section{Fragmentations with infinite reproduction measure}
\label{InFr}
We sketch how the preceding results can be extended to a class
of  self-similar conservative or dissipative fragmentations with infinite reproduction measure.
Such processes were 
introduced in \cite{Be1} where 
the reproduction law was called 
`dislocation measure'.

\par Let
$\nu$ be a measure on the infinite simplex 
$$\Delta=\{(s_j):  s_j\geq 0,\, s_j\downarrow 0, \,\sum_{j=1}^{\infty} s_j\leq 1\}$$ 
such that the integral
\begin{equation}\label{psi-inf}
\psi(\beta):=\int_\Delta (1-\sum s_j^\beta)\,\nu({\rm d}s)
\end{equation}
satisfies $1<\psi(\beta)<\infty$ for some $\beta>0$. 
For $\alpha\geq 0$ we can define a fragmentation process $(X(t), t\geq 0)$ 
 with the property that a generic particle of size $x$ gives birth to
a collection of particles of sizes $xs_j$ with $s=(s_j)\in B$ at rate $x^{\alpha}\nu(B)$,
where $s_j$ runs over nonzero coordinates of $s$ and
$B\subset\Delta$ runs over Borel sets of finite $\nu$-measure.

\par If $\nu$ is a probability measure it can be regarded as a reproduction law by 
defining $(\xi_j)$ to be a random element of $\Delta$ with distribution $\nu$. 
In this case the structural measure is identified
with the superposition of marginal distributions of $\nu$, and the definition (\ref{psi-inf}) agrees with
our definition of the characteristic function in Section \ref{Defin}. 
The case $\nu<\infty$ is easily reduced to the case $\nu(\Delta)=1$ by the obvious time-change. 

\par In the case 
$\nu(\Delta)=\infty$ some features of the fragmentation process are different, in particular, 
each particle produces infinitely many generations within arbitrarily small time period.
As a consequence, the life-time of a particle is not a well-defined quantity, we do not
have a tree representation for the genealogical structure, nor equation like (\ref{d-ident-Mt}).
Still, we can define $\beta_a$ and $\beta^*$ exactly as in the case of finite measure,
and consider $\beta$-sizes for $\Re\,\beta>\beta_a$.
Formula (\ref{pow-ser})  
remains valid and can be proved by
an argument exploiting approximation of $\nu$ by suitable finite measures,
or by using the methods developed in \cite{BC, BY1}.
For $\alpha>0$ conclusions of Theorems \ref{weak-conv} and \ref{strongconv} remain valid 
if we assume that
$$\psi'(\beta^*+)<\infty\quad \hbox{and}\quad 
\int_{\Delta}\left(\sum_{j=1}^\infty s_j^{\beta^*}-1\right)^2\nu({\rm d}s)<\infty\,,$$
and impose a non-arithmeticity condition on $\nu$.
The analog of Corollary \ref{Cint} is shown by arguments similar to those in
Theorem 2 of \cite{Be2}, and the analog of distributional equation  (\ref{FPeqn})
generalises in the form of the Laplace transform identity
\begin{equation}\label{FPeqnInf}
{\cal L}(\theta)=\int_{\Delta} \,\prod_{j=1}^{\infty}\, {\cal L}( \theta s_j^{\beta^*})\,\nu({\rm d}s)\,
\end{equation}
for ${\cal L}(\theta):={\mathbb E}\exp(-\theta M_\infty)$. The question about the uniqueness 
of solution to 
(\ref{FPeqnInf}) with infinite measure seems
to have not been considered before and remains open. 

\par  Finally, we observe the obvious interchangeability between parameters $\alpha$ and $\beta$ (which holds no matter whether the reproduction measure is finite or infinite). In the variables $\widetilde{\xi}_j=\xi^\beta_j$ the fragmentation process 
has the life-time parameter $\alpha/\beta$ and differs only  by a particle-wise transformation of sizes. This observation enables us to handle situations
when the measure $\nu$ is not supported by the infinite simplex $\Delta$, but rather by
$$\Delta_{\beta}=\{(s_j):  s_j\geq 0,\, s_j\downarrow 0, \,\sum_{j=1}^{\infty} s_j^{\beta}\leq 1\}$$ 
for some $\beta>0$. Note also that the change of variables with $\beta<0$ yields mathematically equivalent (though physically curious) process where 
the individual `fragments' grow, but the decaying life-times slow down the total increase of size.


\begin{thebibliography}{99}


\bibitem{A} D. J. Aldous (1999). Deterministic and stochastic models for coalescence
(aggregation, coagulation): a review of the mean-field theory for probabilists. {\it Bernoulli}
{\bf 5}, 3--48.

\bibitem{AP} D. J. Aldous  and J. Pitman  (1998). The standard additive coalescent.
{\it Ann. Probab.} {\bf 26}, 1703-1726.

\bibitem{BG} Y. Baryshnikov and A. Gnedin (2001). Counting intervals in
the packing process. {\it Ann. Appl. Probab. \bf 11}, 863-877.

\bibitem{BK} E. Ben-Naim and P.L. Krapivsky (2000).
Fragmentation with a steady source
{\it Phys. Lett. A \bf 275}, 48-53.

\bibitem{Be} J. Bertoin (2000). A fragmentation process connected to Brownian motion. 
 {\it Probab. Theory Relat. Fields} {\bf 117}, 289-301.


\bibitem{Be0} J. Bertoin (2001).  Homogeneous fragmentation processes, {\it Probab. Theory
Relat. Fields}  {\bf 121}, 301-318.

\bibitem{Be1} J. Bertoin (2002). Self-similar fragmentations. 
{\it Ann. Inst. Henri
Poincar\'e} {\bf 38}, 319-340.

\bibitem{Be2} J. Bertoin (2003). The asymptotic behavior of fragmentation
processes. 
{\it J. Euro. Math. Soc.} {\bf 5}, 395-416.

\bibitem{Be3} J. Bertoin (2004). On small sizes in self-similar fragmentations.
{\it  Stochastic Process. Appl.} {\bf 109}, 13-22.

\bibitem{BC} J. Bertoin and M.-E. Caballero (2002). Entrance from $0+$
for increasing semi-stable Markov processes.  {\it Bernoulli} {\bf 8}, 195-205.


\bibitem{BR} J. Bertoin and A. Rouault (2003). Discretization methods for
homogeneous fragmentations. Preprint.

\bibitem{BY1} J. Bertoin and M. Yor (2001). On subordinators, self-similar Markov processes,
and some factorizations of the exponential law. {\it Elect. Commun. Probab. \bf 6}, 95-106. 
Available at {\tt http://www.math.washington.edu/ejpecp/ecp6contents.html}.


\bibitem{BY} J. Bertoin and M. Yor (2004). On the
exponential functionals of L\'evy processes. In preparation.


\bibitem{BCP} D. Beysens,  X. Campi, and E. Pefferkorn. (1995). {\em Fragmentation
Phenomena}. World Scientific, Singapore.


\bibitem{Bi1} J. D. Biggins (1977). Martingale convergence in the branching random walk. {\it J. Appl. Probability \bf 14}, no. 1, 25--37.

\bibitem{Bi2} J. D. Biggins (1992). Uniform convergence of martingales in
the branching random walk {\it Ann. Probab.} {\bf 20}, 137-151.

\bibitem{BDI} M. D. Brennan and R. Durrett (1986).
Splitting intervals, {\it Ann. Probab.} {\bf 14}, 1024-1036.


\bibitem{BDII} M. D. Brennan and R. Durrett (1987).
Splitting intervals II. Limit laws for lengths. 
{\it Probab. Theory Related Fields \bf 75}, 109--127. 

\bibitem{Filippov} A. F. Filippov (1961). 
On the distribution of the sizes of particles which undergo splitting. {\it Th. Probab. Appl. \bf 6}, 275-293. 

\bibitem{Iserles} A. Iserles and  Y. Liu (1997). Integro-differential equations and
generalized hypergeometric functions, {\it J. Math. Anal. Appl.} {\bf 208}, 404-424.


\bibitem{Jag} P. Jagers (1989). General branching processes as Markov fields.
{\it Stochastic Process. Appl.} {\bf 32}, 183--212.

\bibitem{Kolmogorov} A.N. Kolmogorov (1941). {\"U}ber das logarithmisch normale
Verteilungsgesetz der Dimensionen der Teilchen bei Zerst{\"u}kelung, {\it Soviet Doklady} {\bf
31}, 99-101.

\bibitem{KBNG}
P. L. Krapivsky, E. Ben-Naim and I. Grosse (2004). 
Stable distributions in stochastic fragmentation.
{\it J. Phys. A \bf 37}, 2863-2880.

\bibitem{La} J. Lamperti (1972). Semi-stable Markov processes. {\it Z.
Wahrscheinlichkeitstheorie verw. Gebiete} {\bf 22}, 205--225.

\bibitem{Ly} R. A. Lyons (1997). Simple path to Biggins' martingale convergence
for branching random walk. In: Athreya, K. B.,  and Jagers, P.
(editors).
{\it Classical and Modern Branching Processes.}
Springer-Verlag, New York,  pp. 217-221.

\bibitem{LPP} R. Lyons, R. Pemantle,  and Y. Peres (1995). Conceptual proofs
of {\it L} log {\it L} citeria for mean behaviour of branching
processes. {\it Ann. Probab.} {\bf 23}, 1125-1138.


\bibitem{Marichev} Marichev, O. (1983) {\it Handbook of Integral Transforms of Higher Transcendental
Functions: Theory and Algorithmic Tables}, Ellis Horwood, Chichester.


\bibitem{Mauldin}  R.D. Mauldin  and S.C. Williams (1986) Random recursive constructions:
asymptotic geometric and topological properties, {\it Trans. Amer. Math. Soc.} {\bf 295},
325-426.


\bibitem{Miermont} G. 
Miermont (2004). Self-similar fragmentations derived from the stable tree II: splitting at
hubs. {\it Probab. Theory Relat. Fields} (To appear).

\bibitem{MS} G. Miermont and J. Schweinsberg (2003). Self-similar fragmentations and
stable subordinators. In: {\it S\'eminaire de Probabilit\'es XXXVII}, Lecture Notes in Maths.
1832, pp. 333-359. Springer, Berlin.


\bibitem{Ner} O. Nerman (1981).
On the convergence of supercritical general (C-M-J) branching processes.
{\it Z. Wahrsch. Verw. Gebiete} {\bf 57}, 365--395.

\bibitem{Nev} J. Neveu (1987). Multiplicative martingales for spatial branching
processes.
In {\it Seminar on Stochastic Processes,} 
Progr. Probab. Statist. {\bf 15} pp. 223--242. Birkh\"auser, Boston.

\bibitem{Roesler} U. R{\"o}sler (1992) A fixed point theorem for distributions, {\it Stoch. Proc. Appl.}
{\bf 42}, 195-214.

\bibitem{Ross} S.M.  Ross (1983). {\it Stochastic Processes}, Wiley, N.Y.

\bibitem{Sch} J. Schweinsberg (2001).
Applications of the continuous-time ballot theorem 
to Brownian motion and related processes. {\it
Stochastic Process. Appl. \bf 95}, 151--176.

\bibitem{Uch}  K. Uchiyama  (1982). Spatial growth of a branching process of
particles living in ${\bf R}\sp{d}$. {\it Ann.
Probab. \bf 10}, 896-918.


\end{thebibliography}
\end {document}